# Golomb rulers with real marks: the real and integer optima agree up to ten marks

*Esther Notik · ORCID 0009-0005-4076-6531 · preprint, 25 November 2026*

## Abstract

A Golomb ruler is a set of integer marks whose pairwise differences are distinct. We let the marks be real numbers and ask that every difference be at least 1 and that any two differences be at least 1 apart. Every integer Golomb ruler satisfies these conditions, so the least length $G_R(n)$ of such a real ruler with n marks is at most the optimal Golomb length $G_Z(n)$. Fixing the order of all differences splits the real problem into finitely many polyhedral chambers. A branch-and-bound over these chambers, with an exact rational certificate at every leaf, proves $G_R(n) = G_Z(n)$ for $n \le 10$. Each proof takes under a minute and is verified by an independent checker that uses only exact arithmetic. An exhaustive exact enumeration of all chamber vertices confirms the result for $n \le 7$, and shows that there the shortest real rulers are exactly the optimal integer rulers and their mirror images. For every n we have $n(\sqrt{n} - 1)^2 \le G_R(n) \le G_Z(n)$, essentially a bound of Erdős and of Cilleruelo and Ruzsa, so that $G_R(n) \sim G_Z(n) \sim n^2$. Non-integral chamber vertices first appear at $n = 7$: there are 84, all half-integral. Already there, 80 of them are shortest points of their chambers, 40 of them the only one, although the global optimum is integral; a twelve-mark example of the same kind comes from Singer's construction. On the other hand, every Singer ruler with at most 104 marks, cut at the largest gap, is the only shortest point of its own chamber; we conjectured this for all prime powers, and Ma and Zhang have since proved it. The programs, the proofs and the checkers are published with this note at doi:10.5281/zenodo.23843508.

## 1. Introduction

A *Golomb ruler* with n marks is a set of integers $0 = x_1 < x_2 < \dots < x_n$ whose differences $x_j - x_i$ $(i < j)$ are pairwise distinct. Let $G_Z(n)$ be the least possible length $x_n$. For n = 2, …, 10 it is 1, 3, 6, 11, 17, 25, 34, 44, 55 (OEIS A003022 [13]). It is known for $n \le 28$, in each case by exhaustive computer search; the last cases were completed by distributed.net in 2022 [5].

**Definition 1.** A *real 1-separated ruler* with n marks is a sequence of real numbers $0 = x_1 < \dots < x_n$ with $x_j - x_i \ge 1$ for all $i < j$ and $|(x_j - x_i) - (x_l - x_k)| \ge 1$ for all pairs $(i, j) \ne (k, l)$ with $i < j$ and $k < l$. Let $G_R(n)$ be the infimum of $x_n$ over such rulers.

For integer marks, distinct differences are automatically at least 1 apart, so every Golomb ruler is a real 1-separated ruler and $G_R(n) \le G_Z(n)$. The question is whether real marks can ever do better. The same relaxation of the Erdős distinct subset sums problem was studied in [12], where the real and integer optima agree for $n \le 8$. For rulers the structure is simpler, because every condition involves only two differences.

**Earlier work.** For $n \ge 3$ the condition that the differences are at least 1 follows from the other one, so $G_R(n)$ is the least diameter of n points on a line whose pairwise distances differ by at least 1. Erdős [24, Problem 20] conjectured that n points of a Euclidean space with this property have diameter at least $(1 + o(1))n^2$, and proved this on the line; see Erdős Problem #670 [25] and Ho [26], whose quantity $\Delta_1(n)$ is $G_R(n)$. Cilleruelo and Ruzsa [18] studied the problem on the line in 2004 in a scaled form. For n real numbers in [0, 1] they consider the least distance between pair sums, whose largest possible value is $1/G_R(n)$, and prove $G_R(n) \ge n(n - 2\sqrt{n})$. The author learned of these sources only after the earlier versions of this note, which did not mention them (see Sections 8 and 9).

**Results.** Real and integer optima agree for every $n \le 10$ (Theorem 1). The proof is a certified branch-and-bound (Section 3). For $n \le 7$ an exhaustive enumeration of all chamber vertices gives an independent proof, identifies all shortest real rulers (Theorem 2), and shows the first non-integral vertices. Section 6 shows that most of these are shortest points of their chambers, with a short certificate for one of them (Proposition 1), and gives a twelve-mark example from Singer's construction (Proposition 2). So chamber optima can be non-integral even where the global optimum is integral, and Section 6 also shows that no universal floor/ceiling rounding explains the integral optima. Section 7 shows that every Singer ruler for a prime power up to 103 is the only shortest point of its chamber, and conjectures this in general (since proved, see the note in Section 10). Section 8 proves $n(\sqrt{n} - 1)^2 \le G_R(n)$ for all n (Theorem 3), essentially the bound of [24] and [18], so that $G_R(n) \sim G_Z(n) \sim n^2$.

## 2. Chambers, vertices and valid inequalities

Write $k = n - 1$ and use the gaps $g_t = x_{t+1} - x_t$ as coordinates. Every difference $x_j - x_i$ equals $1_I \cdot g$ for the interval $I = \{i, \dots, j - 1\}$, where $1_I$ is its indicator vector. Hence every condition of Definition 1 reads $|e \cdot g| \ge 1$ for e in the *relation set* $R_n$: the interval indicators $1_I$ (among them the unit vectors, giving $g_t \ge 1$) and the differences $1_I - 1_J$ of two distinct intervals, each normalized so that its first nonzero entry is +1, with repetitions removed. $|R_n| = 4, 12, 30, 65, 126$ for $n = 3, \dots, 7$.

For a sign vector σ on $R_n$ with $\sigma_e = +1$ on the unit vectors, the *chamber* $P_\sigma$ is the set of g with $\sigma_e\, e \cdot g \ge 1$ for all $e \in R_n$. The real 1-separated rulers form the disjoint union of the nonempty chambers. A chamber corresponds to a choice of the order of all differences, which is a combinatorial type of real Golomb ruler in the sense of Beck, Bogart and Pham [2].

**Lemma 1 (vertices).** $G_R(n)$ is attained. It is the least total length $g_1 + \dots + g_k$ over the finite set of points g such that (i) $g = M^{-1}s$ for a $k \times k$ matrix M whose rows are linearly independent elements of $R_n$ and a sign vector $s \in \{-1, 1\}^k$, and (ii) $g_t \ge 1$ for all t and $|e \cdot g| \ge 1$ for all $e \in R_n$.

*Proof.* A nonempty chamber is a polyhedron contained in $\{g : g \ge 1\}$. It contains no line, and the length is bounded below on it and strictly increasing along every recession direction. By the fundamental theorem of linear programming its minimum over the chamber is attained, at a vertex, where k linearly independent defining inequalities hold with equality: $e \cdot g = \sigma_e$ for k independent rows e. Such a point has the form (i), and it satisfies (ii) because it lies in the chamber. Conversely, every point satisfying (ii) is a real 1-separated ruler. ■

By Cramer's rule $|\det M| \cdot g$ is an integer vector. A feasible vertex is *integral* if g is an integer vector; it is then an integer Golomb ruler.

**Lemma 2 (packing).** Any r distinct differences of a real 1-separated ruler have sum at least $1 + 2 + \dots + r = r(r + 1)/2$, because they are at least 1 and pairwise at least 1 apart. In particular the largest difference, $x_n$, is at least $n(n - 1)/2$. Since {0, 1}, {0, 1, 3} and {0, 1, 4, 6} are Golomb rulers, $G_R(n) = G_Z(n) = n(n - 1)/2$ for $n \le 4$.

**Lemma 3 (sub-rulers).** For $i < j$ the marks $x_i, \dots, x_j$ form, after translation, a real 1-separated ruler with $m = j - i + 1$ marks, so $x_j - x_i \ge G_R(m)$.

**Lemma 4 (symmetry).** The mirror image $x \mapsto x_n - x$ of a real 1-separated ruler is again one, of the same length, with the gaps in reverse order. The gaps $g_1$ and $g_k$ ($k \ge 2$) are two distinct differences, so $|g_1 - g_k| \ge 1$. Hence a lower bound proved for all rulers with $g_k - g_1 \ge 1$ holds for all rulers.

## 3. The certified branch-and-bound

To prove $G_R(n) \ge T$ we search a binary tree. A node is described by its path, a list of branching rows $\sigma\, e \cdot g \ge 1$ with $e \in R_n$; the root has an empty path. At each node we solve the linear program: minimize $g_1 + \dots + g_k$ subject to the rows valid at the node, which are

- **$u_t$:** $g_t \geq 1$ (Definition 1);
- **$s_I$:** $1_I \cdot g \geq L(m)$ for every interval I of m marks with $3 \leq m \leq n - 1$ (Lemma 3), where $L(m) = m(m - 1)/2$ for $m \leq 4$ (Lemma 2) and $L(m) = G_Z(m)$ for $5 \leq m \leq n - 1$, as certified by the proofs for smaller n;
- **p:** packing cuts $\Sigma_{I \in S} 1_I \cdot g \geq r(r + 1)/2$ for a set S of r distinct intervals (Lemma 2), added when the LP point violates one;
- **y:** $g_k - g_1 \geq 1$ (Lemma 4);
- **$b_d$:** the branching rows of the path.

If the LP value is at least T, the node is a leaf and receives a certificate: multipliers $y_i \geq 0$ with $\Sigma\, y_i\, a_i = (1, \ldots, 1)$ and $\Sigma\, y_i\, b_i \geq T$ for rows $a_i \cdot g \geq b_i$ valid at the node. For every point of the node, the length then equals $\Sigma\, y_i\, a_i \cdot g \geq \Sigma\, y_i\, b_i \geq T$. (An empty node would instead receive a Farkas certificate, $\Sigma\, y_i\, a_i = 0$ with $\Sigma\, y_i\, b_i > 0$; none was needed.) Otherwise we add a violated packing cut if there is one, and then branch on the relation $e \in R_n$ with the smallest $|e \cdot g|$ at the LP point: the children add $e \cdot g \geq 1$ and $-e \cdot g \geq 1$. Every ruler satisfies $|e \cdot g| \geq 1$, so the two children cover the node. If the LP point satisfied every relation, it would be a real ruler shorter than T, and the proof would fail.

The LP is solved in floating point by HiGHS 1.15.1 [17]. The certificates are exact: the k rows that are tight in HiGHS's optimal basis are taken, and the multipliers are obtained from them by solving a k × k linear system in rational arithmetic. They are kept only if they are nonnegative and give $\Sigma\, y_i\, b_i \geq T$ exactly. Floating point therefore only guides the search; the proof rests on the certificates.

**The checker** (Real_Golomb_rulers_check.py) uses only the Python standard library and exact fractions. For a proof file it recomputes $R_n$ and the intervals from n and compares them with the file. It checks that the tree has one root, that every internal node has exactly the two children $e \cdot g \geq 1$ and $-e \cdot g \geq 1$ for a relation $e \in R_n$, and that every leaf has a valid certificate. It checks that every row used is one of the kinds listed above, with an index in range, and with each packing cut and sub-ruler bound within what Lemmas 2 and 3 allow. Sub-ruler bounds for $m \geq 5$ are accepted only together with the proof for m, which is checked recursively. Finally, it checks an optimal integer ruler for the upper bound. Checking the proof for n = 10, with its whole chain for n = 5, …, 9, takes about 3 seconds.

**Tests.** With the target raised by 1 ($T = G_Z(n) + 1$), the search must fail. For n = 5, …, 8 it does: it finds exactly the optimal integer rulers, one of each mirror pair. The checker rejects a proof with open leaves, a proof in which one multiplier was changed, and proofs with malformed or out-of-range row references (for example a row $g_t \geq 1$ with t outside 1, …, k, which would read $0 \geq 1$). For $n \leq 7$ the result agrees with the enumeration of Section 5.

**Table 1.** The certified proofs of $G_R(n) \geq G_Z(n)$, each run on one process of a laptop (Intel Core Ultra 7 155U).

| n | target | nodes | leaves | packing cuts | max depth | time | proof file |
|---|---|---|---|---|---|---|---|
| 5 | 11 | 11 | 6 | 7 | 4 | 0.01 s | 2.5 KB |
| 6 | 17 | 47 | 24 | 16 | 6 | 0.03 s | 8.7 KB |
| 7 | 25 | 643 | 322 | 55 | 15 | 0.4 s | 94 KB |
| 8 | 34 | 3,823 | 1,912 | 128 | 23 | 3.3 s | 0.57 MB |
| 9 | 44 | 15,929 | 7,965 | 298 | 32 | 19 s | 2.5 MB |
| 10 | 55 | 18,425 | 9,213 | 332 | 29 | 33 s | 3.1 MB |

## 4. Results

**Theorem 1.** For $2 \leq n \leq 10$, $G_R(n) = G_Z(n)$, that is, 1, 3, 6, 11, 17, 25, 34, 44, 55.

*Proof.* For $n \le 4$ this is Lemma 2. For $5 \le n \le 10$ the proof files of Table 1 give $G_R(n) \ge G_Z(n)$. They are built in the order n = 5, 6, …, 10, each using the bounds for smaller m as sub-ruler rows. The upper bound comes from the optimal integer rulers {0, 1, 4, 9, 11}, {0, 1, 4, 10, 12, 17}, {0, 1, 4, 10, 18, 23, 25}, {0, 1, 4, 9, 15, 22, 32, 34}, {0, 1, 5, 12, 25, 27, 35, 41, 44} and {0, 1, 6, 10, 23, 26, 34, 41, 53, 55}. ■

The proof is computer-assisted. Every step can be checked with the published checker, without the LP solver.

## 5. Exhaustive enumeration of the vertices for n ≤ 7

Lemma 1 also gives a direct method. Enumerate every k-element subset of $R_n$. For each subset, fraction-free Gauss–Jordan elimination in exact integer arithmetic gives $D = \pm\det M$ and $A = D{\cdot}M^{-1}$; every division is checked to be exact, and $M{\cdot}A = D{\cdot}I$ is checked afterwards. All integers involved are small: in Bareiss's fraction-free elimination [1] every intermediate entry is, up to sign, a minor of [M | I], hence by Hadamard's inequality at most $6^3 = 216$ in absolute value for n = 7, so no numerator exceeds $2{\cdot}216^2 < 10^5$ and the 64-bit arithmetic cannot overflow; the floating-point product A s has integer partial sums below 400, so it is exact. For each sign vector s, keep the point g = A s / D if it satisfies (ii). We ran this in pure Python for $n \le 6$ and in a vectorized numpy version for n = 7: 4,925,156,775 subsets, 6.3 hours on 10 worker processes. The numpy version reproduces every count of the Python version for n = 5 and 6.

Every feasible vertex lies in exactly one chamber, and every chamber contains a feasible vertex (Lemma 1). The vertices found meet 2, 10, 114, 2,608 and 107,498 chambers for n = 3, …, 7. These are exactly the numbers of combinatorially different real Golomb rulers listed by Beck, Bogart and Pham [2, Section 5], which also form OEIS A237749 [13]. So the enumeration met every chamber.

**Table 2.** The enumeration. "Solutions" are the (subset, sign) pairs that give a feasible vertex, by |det M|. "Chambers" is the number of distinct chambers containing a feasible vertex.

| n | $\lvert R_n \rvert$ | subsets | nonsingular | solutions by \|det M\| | feasible vertices | chambers | non-integral |
|---|---|---|---|---|---|---|---|
| 3 | 4 | 6 | 6 | 1: 2 | 2 | 2 | 0 |
| 4 | 12 | 220 | 189 | 1: 10 | 10 | 10 | 0 |
| 5 | 30 | 27,405 | 21,035 | 1: 176 | 118 | 114 | 0 |
| 6 | 65 | 8,259,888 | 5,903,710 | 1: 6,296; 2: 2 | 3,144 | 2,608 | 0 |
| 7 | 126 | 4,925,156,775 | 3,365,649,124 | 1: 425,802; 2: 1,378; 3: 2 | 165,522 | 107,498 | 84 |

**Theorem 2.** For n = 5, 6, 7 the real 1-separated rulers of length $G_Z(n)$ are exactly the optimal Golomb rulers below and their mirror images:

- n = 5: {0, 1, 4, 9, 11} and {0, 2, 7, 8, 11};
- n = 6: {0, 1, 4, 10, 12, 17}, {0, 1, 4, 10, 15, 17}, {0, 1, 8, 11, 13, 17} and {0, 1, 8, 12, 14, 17};
- n = 7: {0, 1, 4, 10, 18, 23, 25}, {0, 1, 7, 11, 20, 23, 25}, {0, 1, 11, 16, 19, 23, 25}, {0, 2, 3, 10, 16, 21, 25} and {0, 2, 7, 13, 21, 22, 25}.

*Proof.* The enumeration finds the least length over the finite set of Lemma 1, attained only at these rulers and their mirror images (4, 8 and 10 vertices). They lie in pairwise different chambers. In each chamber the length is minimized on a bounded face, the convex hull of the optimal vertices of that chamber, so each of these chambers has a unique shortest ruler, and no other real ruler has length $G_Z(n)$. The next shortest feasible vertices have lengths 12, 18 and 26. ■

# 6. Non-integral vertices and optima

Up to n = 6 every feasible vertex is integral. At n = 6, two solutions come from subsets with |det M| = 2, but their vertices are integral and are also reached by subsets with |det M| = 1. At n = 7 the first non-integral feasible vertices appear: 84 of them, all half-integral. Two solutions come from subsets with |det M| = 3, but no vertex has denominator 3. The shortest non-integral vertex has length 33, far above $G_Z(7)$ = 25:

0, 7/2, 6, 11, 20, 21, 33, with differences 1, 5/2, 7/2, 5, 6, 15/2, 9, 10, 11, 12, 13, 14, 15, 33/2, 35/2, 20, 21, 22, 27, 59/2, 33, which are pairwise at least 1 apart.

So the simple route to $G_R(n) = G_Z(n)$ for all n, "every vertex is integral", fails. A sharper question is whether a non-integral vertex can be the shortest point of its chamber. For each of the 84 vertices we solved the linear program of its chamber and recorded an exact certificate. For 80 of them a multiplier vector $y \geq 0$ with $\Sigma\, y_i a_i = (1, \ldots, 1)$ over rows of the chamber has $\Sigma\, y_i$ equal to the length of the vertex, so the vertex is a shortest point of its chamber. For 40 of these the rows with $y_i > 0$ have rank 6, so the vertex is the only shortest point; for the other 40 a second point of the chamber with the same length is given. Each of the remaining 4 comes with a shorter point of its chamber (Table 3). These verdicts do not depend on the completeness of the enumeration; only the list of 84 vertices does.

**Table 3.** The chambers of the 84 non-integral vertices at n = 7.

| the vertex is | vertices | shortest of them |
|---|---|---|
| the only shortest point of its chamber | 40 | 42 |
| a shortest point, not the only one | 40 | 33 |
| not a shortest point | 4 | 38 |

So non-integral chamber optima occur already at n = 7, the first n with non-integral vertices, although $G_R(7)$ = 25 is attained only at integer rulers (Theorem 2). The shortest of them, 0, 7/2, 6, 11, 20, 21, 33, shares the optimal face of its chamber with the integer ruler {0, 4, 6, 11, 20, 21, 33}. The same happens for separated subset sums [12]: some chambers there have non-integral optimal vertices, but none beats the integer optimum for n ≤ 8. The shortest unique example has a short certificate.

**Proposition 1 (seven marks).** The real 1-separated ruler P = 0, 9, 14, 33/2, 20, 32, 42 is the only shortest point of its chamber.

*Proof.* The 21 differences of P are 5/2, 7/2, 5, 6, 15/2, 9, 10, 11, 12, 14, 31/2, 33/2, 18, 20, 22, 23, 51/2, 28, 32, 33, 42. They are at least 1 and pairwise at least 1 apart, so P lies in its chamber. Table 4 lists six rows of the chamber that are tight at P, with multipliers 8, 4, 7, 6, 9, 8. Their weighted sum is $g_1 + \ldots + g_6 \geq 42$, so no point of the chamber is shorter than 42, and a point of length 42 satisfies the six rows with equality. The six rows have determinant 2, so that point is P. ■

**Table 4.** The certificate of Proposition 1. Each row says that two differences that are consecutive in the order of P differ by at least 1; the last column gives their values at P.

| row of the chamber | multiplier | differences at P |
|---|---|---|
| $-g_3 + g_4 \geq 1$ | 8 | 5/2 → 7/2 |
| $-g_2 + g_3 + g_4 \geq 1$ | 4 | 5 → 6 |
| $-g_1 + g_6 \geq 1$ | 7 | 9 → 10 |
| $g_2 + g_3 + g_4 - g_6 \geq 1$ | 6 | 10 → 11 |
| $-g_2 - g_3 - g_4 + g_5 \geq 1$ | 9 | 11 → 12 |
| $g_1 + g_2 + g_3 - g_4 - g_5 \geq 1$ | 8 | 31/2 → 33/2 |

With twelve marks there is an example from Singer's construction.

**Proposition 2 (twelve marks).** Let X = {0, 4, 24, 27, 67, 77, 86, 88, 103, 119, 125, 132}, a Golomb ruler with n = 12 marks. The only shortest point of its chamber is the half-integral ruler

P = 0, 4, 22, 25, 123/2, 141/2, 157/2, 161/2, 94, 217/2, 114, 121, of length 121.

*Proof.* Let $d_{(1)} < d_{(2)} < \dots < d_{(66)}$ be the differences in the order they have in X. The chamber of X is the set of g with $d_{(1)}(g) \geq 1$ and $d_{(r+1)}(g) - d_{(r)}(g) \geq 1$ for r = 1, …, 65; these rows imply all others. The differences of P come in the same order as those of X and are pairwise at least 1 apart, so P lies in the chamber. Exactly 24 rows are tight at P. Table 5 lists them with multipliers $y_i \geq 1/2$, for which $\Sigma\, y_i\, a_i = (1, \dots, 1)$ and $\Sigma\, y_i = 121$ in exact arithmetic. Hence every point g of the chamber has length $\Sigma\, y_i\, a_i \cdot g \geq \Sigma\, y_i = 121$, with equality only if all 24 rows are tight at g. These rows have rank 11, so g = P. ■

**Table 5.** The 24 rows tight at P in Proposition 2. Each row says that two differences that are consecutive in the order of X differ by at least 1. It is given by the two differences in X and in P, and by its multiplier y.

| in X | in P | y | in X | in P | y |
|---|---|---|---|---|---|
| 2 → 3 | 2 → 3 | 25/2 | 23 → 24 | 21 → 22 | 1/2 |
| 3 → 4 | 3 → 4 | 1/2 | 29 → 31 | 27 → 28 | 1/2 |
| 7 → 9 | 7 → 8 | 19/2 | 36 → 37 | 65/2 → 67/2 | 1/2 |
| 9 → 10 | 8 → 9 | 11/2 | 39 → 40 | 71/2 → 73/2 | 25/2 |
| 10 → 11 | 9 → 10 | 1/2 | 43 → 44 | 79/2 → 81/2 | 23/2 |
| 13 → 15 | 25/2 → 27/2 | 1/2 | 46 → 48 | 85/2 → 87/2 | 1/2 |
| 15 → 16 | 27/2 → 29/2 | 57/2 | 58 → 59 | 105/2 → 107/2 | 1/2 |
| 16 → 17 | 29/2 → 31/2 | 59/2 | 61 → 62 | 111/2 → 113/2 | 1/2 |
| 19 → 20 | 17 → 18 | 1/2 | 62 → 63 | 113/2 → 115/2 | 1/2 |
| 20 → 21 | 18 → 19 | 1/2 | 63 → 64 | 115/2 → 117/2 | 1/2 |
| 21 → 22 | 19 → 20 | 7/2 | 64 → 65 | 117/2 → 119/2 | 1/2 |
| 22 → 23 | 20 → 21 | 1/2 | 98 → 99 | 89 → 90 | 1/2 |

X is a cut of a Singer perfect difference set modulo 133 [16], cut at a gap that is not the largest; it was found while testing whether Singer rulers are the shortest points of their chambers. Its mirror image gives the mirror case. Both examples are far above the integer optimum (42 against $G_Z(7) = 25$, and 121 against $G_Z(12) = 85$), so they do not bear on Theorem 1.

**Remark (rounding).** Simple rounding does not account for the integral global optima. Rounding each of the five half-integral marks of P in Proposition 2 to its floor or its ceiling never gives a Golomb ruler of length at most 121. At n = 14, the half-integral real 1-separated ruler

0, 3, 55/2, 48, 57, 58, 90, 104, 279/2, 283/2, 157, 163, 168, 176,

a shortest point of the chamber of another Singer cut, is not even within 1 of an integer Golomb ruler of length at most 176: no such ruler has every mark between the floor of the corresponding mark minus 1 and its ceiling plus 1. Among 1,936 non-integral chamber optima and vertices for n ≤ 14 (the 84 vertices of Table 3, four minima from the Singer test of Section 7, and chamber optima found by a walk through neighbouring chambers), rounding to floors and ceilings fails for 1,050. These examples rule out a universal length-nonincreasing floor/ceiling rounding procedure, even one that may move each mark one step further.

# 7. Singer rulers

For a prime power m, Singer [16] constructs a set D of m + 1 residues modulo $q = m^2 + m + 1$ whose differences cover every non-zero residue exactly once, a perfect difference set; so does every multiple t·D with gcd(t, q) = 1. Cutting the circle $\mathbb{Z}/q\mathbb{Z}$ at a gap between two consecutive elements of t·D and reading the elements as integers in [0, q) gives a Golomb ruler with n = m + 1 marks, shortest when the cut is at the largest gap. We call these rulers, cut at the largest gap, *Singer rulers*. The shortest Singer ruler for each m has length $G_Z(n)$ for n = 3, 4, 5, 6, 10, 12, 14, 18, 20, 24, 26 and 28, that is, for every n ≤ 28 with n − 1 a prime power except n = 8, 9 and 17, where it is longer by 1, 1 and 2. The ruler of Theorem 1 for n = 10 is itself a perfect difference set modulo 91.

**Computation.** For every prime power m ≤ 103 (37 values of m, up to 104 marks) and every multiplier t, we solved the linear program of the chamber of the Singer ruler and certified the result exactly as in Section 6: in each case a multiplier vector y with $\Sigma\, y_i a_i = (1, \ldots, 1)$ and $\Sigma\, y_i$ equal to the length, positive on rows of rank n − 1. So each of these 24,450 distinct rulers is the only shortest point of its own chamber, including those for n = 8, 9 and 17, which are not globally optimal. For m ≤ 47 we also tested the cuts at all other gaps (116,868 rulers in all), again with every verdict certified exactly: from m = 19 on each of them is the only shortest point of its chamber, while for m ≤ 17 there are 396 that are not shortest points and 18 that are shortest but not the only one, each with an exact second point of the same length. Proposition 2 comes from one of these.

**Conjecture.** For every prime power m, every Singer ruler is the only shortest point of its chamber.

Chamber optimality is local and does not bound $G_R(n)$ by itself. For the analogous relaxation of distinct subset sums, the Conway–Guy and Bohman sets have structured dual certificates valid for all n [12]; a certificate of that kind for Singer rulers, built from the cyclic structure of D, would prove the conjecture.

# 8. Bounds for all n

**Theorem 3.** For every n ≥ 2, $n(\sqrt{n} - 1)^2 = n^2 - 2n^{3/2} + n \le G_R(n) \le G_Z(n) \le n^2 + o(n^2)$. In particular $G_R(n) \sim G_Z(n) \sim n^2$.

*Proof.* Let $0 = x_1 < \ldots < x_n = L$ be a real 1-separated ruler, fix u > 0, and for real t let k(t) be the number of marks in [t, t + u). Then k(t) = 0 outside an interval of length L + u, and $\int k(t)\,dt = nu$. Two marks at distance d both lie in [t, t + u) for a set of t of length max(0, u − d), so $\int k(t)(k(t) - 1)/2\,dt$ is the sum of max(0, u − d) over the n(n − 1)/2 differences d. These are at least 1 and pairwise at least 1 apart, so the r-th smallest is at least r, and the sum is at most $\Sigma_{r \ge 1} \max(0, u - r) \le u^2/2$. By the Cauchy–Schwarz inequality $\int k(t)^2\,dt \ge (nu)^2/(L + u)$. Hence $n^2u^2/(L + u) - nu \le u^2$, that is, $L \ge n^2u/(n + u) - u$, and the choice $u = n^{3/2} - n$ gives $L \ge n^2 - 2n^{3/2} + n$. For the upper bound let p be the least prime with p ≥ n − 1. A Singer set modulo $p^2 + p + 1$ has p + 1 ≥ n elements whose differences are distinct even modulo $p^2 + p + 1$, so any n of them, read as integers in $[0, p^2 + p]$, form a Golomb ruler of length at most $p^2 + p$. By the prime number theorem p = n + o(n). ■

The lower bound is the counting argument of Erdős and Turán [7] for Sidon sets. It uses only that the differences are at least 1 apart, so it applies to real rulers unchanged; no bound of this counting kind can tell $G_R$ from $G_Z$. Theorem 3 gives $G_Z(n) - G_R(n) = o(n^2)$; whether the difference is always 0 is the main open question.

For real numbers the lower bound is essentially due to Erdős [24, Problem 20]. Summing the differences $x_{i+k} - x_i$ over i and over $k \le \sqrt{n}$, he obtains, roughly, $G_R(n) > (n - \sqrt{n})^2$, and concludes that $G_R(n) \ge (1 + o(1))n^2$. The bound of Theorem 3 is the same expression, since $n(\sqrt{n} - 1)^2 = (n - \sqrt{n})^2$. A bound of the same kind was proved by Cilleruelo and Ruzsa [18, Theorem 1], with the method of Lindström: their statement is $G_R(n) \ge n(n - 2\sqrt{n})$ in our notation. They ask whether $2\sqrt{n}$ can be replaced by $o(\sqrt{n})$, and note that any improvement of the coefficient 2 would improve the bound then known for integer Sidon sets. So the lower bound of Theorem 3 is not new.

# 9. Related work

Beck, Bogart and Pham [2] introduced real Golomb rulers and their combinatorial types in order to count integer rulers with Ehrhart theory for inside-out polytopes; Bogart and Cuéllar [3] extended this to generalized rulers. Neither paper studies the separated minimum $G_R(n)$. Kocuk and van Hoeve [8] compare integer, constraint and quadratic programming methods for certifying optimal integer rulers, up to 14 marks. They make one of their models exact by branching on $d_{ij} \le d_{kl} - 1$ or $d_{ij} \ge d_{kl} + 1$, the same dichotomy as our branching, but inside an integer model with integer bounds. Their bounds $d_{ij} \ge G_{j-i+1}$ are the integer counterpart of our sub-ruler rows. The linear-programming lower bounds of Lorentzen and Nilsen [9] and others, shown to coincide by Meyer and Jaumard [10], rest on packing inequalities like those of Lemma 2, which also hold for real 1-separated rulers. Duxbury, Lavor and Salles-Neto [6] propose a different continuous model: the marks stay at integer positions, and the 0–1 variables that say which positions carry a mark are relaxed to [0, 1]. They conjectured that this non-convex model is exact, on evidence from a local solver; Ouzia [14] presents a proof of the conjecture and corrects the outputs of their greedy heuristic. In that model the marks remain at integer positions, so the question does not bear on $G_R(n)$. Nathanson [11] bounds the size of integer Sidon sets whose sumset is Δ-separated, for a positive integer Δ. The conditions of Definition 1 say exactly that the sums $x_i + x_j$ $(i \le j)$ are pairwise at least 1 apart, so such a set with least element 0, divided by Δ, is a real 1-separated ruler with marks in $(1/\Delta)\mathbb{Z}$, and conversely. Let $D_\Delta(n)$ be the least largest element of such a set with n elements. Since $G_R(n)$ is attained at a rational point (Lemma 1), $G_R(n)$ is the minimum of $D_\Delta(n)/\Delta$ over the positive integers Δ, and $G_R(n) = G_Z(n)$ says exactly that $D_\Delta(n) = \Delta \cdot G_Z(n)$ for every positive integer Δ.

Cilleruelo and Ruzsa [18] studied $G_R(n)$ in the form $1/G_R(n)$, the largest separation of the pair sums of n numbers in [0, 1], together with the analogous problems modulo one, for infinite sequences and for p-adic integers. Erdős [24] and they give the lower bounds quoted in Section 8. They do not determine values of $G_R(n)$ or compare it with $G_Z(n)$ for given n. Erdős's problem [24, 25] concerns the asymptotic size of the diameter in every dimension. In a comment of 28 September 2026 at [25], veljjanoski, comparing planar sets with Golomb rulers, notes not knowing whether real points on a line can be shorter than these. By Theorem 1 they cannot for $n \le 10$. So the planar sets of 5 to 10 points reported there, for example five points of diameter less than 10.187, are smaller than any such set on a line. Versions of this note before 25 November 2026 stated that no earlier study of $G_R(n)$ had been found; that statement is corrected here. We know of no determination of exact values of $G_R(n)$ before this note.

# 10. Open problems

- Is $G_R(n) = G_Z(n)$ for every n? Theorem 3 gives $G_R(n) \sim G_Z(n)$. A first n with $G_R(n) < G_Z(n)$ would need a non-integral chamber optimum below $G_Z(n)$, and by the remark of Section 6 no universal floor/ceiling rounding can rule this out.
- n = 11 and 12. Random-dive estimates of the tree size suggest about $3\text{–}7 \cdot 10^6$ nodes for n = 11 and $0.8\text{–}1.4 \cdot 10^7$ for n = 12, that is, hours on one process; we have not run them.
- Uniqueness for n = 8, 9, 10: are the shortest real rulers again exactly the optimal integer rulers? The proofs of Theorem 1 give only the length.
- Non-integral chamber optima exist from n = 7 on (Section 6), yet no global optimum found so far is non-integral. What keeps them away from the optimum? How close to $G_Z(n)$ can a non-integral chamber optimum come? At n = 7 the shortest has length 33, against $G_Z(7) = 25$.
- Prove the conjecture of Section 7 for all prime powers, for instance by structured dual certificates as for the Conway–Guy and Bohman sets in [12]; and ask the same for the rulers of Bose [4] and Ruzsa [15].

**Note added, 25 November 2026.** The questions above are left as they stood in the version of 17 November. Since then, the following answers have appeared.

- n = 11 and 12. Genkin [19] proves $G_R(11) = G_Z(11) = 72$ and $G_R(12) = G_Z(12) = 85$ by the method of this note, with stronger packing cuts. The certificates were checked for this note with a separately written exact checker, which accepts them. A computation for n = 11 had been posted on 20 November [20], and $G_R(13) = G_Z(13) = 106$ is claimed in [21]; these two computations have not been checked by the author.
- Uniqueness for n = 8, 9, 10. A Lean formalization in the repository of the Omega Institute [22], built on the certificates of this note, reports that the shortest real rulers are exactly the optimal integer rulers and their mirror images. It has not been compiled by the author.
- The conjecture of Section 7 and the last question above. Ma and Zhang [23] prove that every cyclic perfect difference set, cut at its largest gap, gives a ruler that is the only shortest point of its chamber, with explicit dual certificates. This proves the conjecture. For the rulers of Bose and Ruzsa they give counterexamples to the corresponding statement. The proof and the counterexamples were checked for this note. A Lean version of the main step was merged earlier, on 18 November, in the repository of [22] (pull request #31500).

# Data and code

The Zenodo record of this note, doi:10.5281/zenodo.23843508 (https://doi.org/10.5281/zenodo.23843508, which always resolves to the latest version), contains the archive Real_Golomb_rulers_code_data.zip, with these parts and a list of SHA-256 checksums:

- the branch-and-bound (Real_Golomb_rulers_bnb.py, needs numpy and highspy);
- the checker (Real_Golomb_rulers_check.py);
- the six proof files for n = 5, …, 10 (data/bnb/, JSON lines);
- the exact check of Propositions 1 and 2 and of the remark of Section 6 (Real_Golomb_rulers_propositions_check.py), which rebuilds the chambers, checks every entry of Tables 4 and 5, and searches all the roundings of the remark;
- the certificates behind Table 3 (data/n7_nonintegral_chambers.json) and their checker (Real_Golomb_rulers_n7_chambers_check.py), which also recomputes the 84 vertices from the enumeration output;
- the Singer chamber test of Section 7 (Real_Golomb_rulers_singer_chambers.py, needs numpy and highspy) and the rounding search behind the remark of Section 6 (Real_Golomb_rulers_rounding.py), with their outputs (data/singer_chambers/, data/rounding/);
- the exact certificates for the 258 Singer chamber verdicts whose uniqueness the first run decided only in floating point (data/singer_chambers/recheck_2026-11-17.json), the program that produced them (Real_Golomb_rulers_singer_recheck.py, needs numpy and highspy) and their checker (Real_Golomb_rulers_singer_recheck_check.py);
- the vertex enumeration in pure Python and numpy, with a C version that was not used;
- all enumeration outputs (data/real_golomb_runs/).

To check Theorem 1 for n = 10: python -I Real_Golomb_rulers_check.py ../data/bnb/proof_n10_T55.jsonl with --dep m:../data/bnb/proof_n<m>_T<$G_Z$(m)>.jsonl for m = 5, …, 9 and --upper 0,1,6,10,23,26,34,41,53,55. To check Propositions 1 and 2, the remark of Section 6 and Table 3: python -I Real_Golomb_rulers_propositions_check.py and python -I Real_Golomb_rulers_n7_chambers_check.py; for the certificates of Section 7: python -I Real_Golomb_rulers_singer_recheck_check.py. The four checkers use only the Python standard library, and every check raises an error explicitly, so they also check under python -O. The archive README gives the exact commands.

## Acknowledgements

The research was directed by the author, who used Claude (Anthropic) for the programs and a draft of the paper; the author checked the results.